\documentclass[11pt]{article}
\usepackage{amsfonts,amsmath,amssymb,accents,amsthm}%,mathtools}
\usepackage{verbatim}
\usepackage{hyperref}
\newcommand{\dis}{\displaystyle}

\newcommand{\noi}{\noindent}
\newcommand{\halmos}{\rule{1ex}{1.4ex}}
\newcommand{\QED}{\nopagebreak{\hspace*{\fill}$\halmos$\medskip}}

\newcommand{\quand}{\quad\mbox{and}\quad}

\newtheoremstyle{mythm}% name
  {}%      Space above
  {}%      Space below
  {\itshape}%         Body font
  {}%         Indent amount (empty = no indent, \parindent = para indent)
  {\bfseries}% Thm head font
  {}%        Punctuation after thm head
  {.5em}%     Space after thm head: " " = normal interword space;
  {#1 #2 \thmnote{(#3)}}%  Thm head spec (can be left empty, meaning `normal')

\theoremstyle{mythm}
\newtheorem{theorem}{Theorem}%[section]
\newtheorem{proposition}[theorem]{Proposition}
\newtheorem{lemma}[theorem]{Lemma}
\newtheorem{exercise}[theorem]{Exercise}
\newtheorem{corollary}[theorem]{Corollary}
\newtheorem{conjecture}[theorem]{Conjecture}

\newtheorem{counterex}[theorem]{Counterexample}

\newcommand{\bt}{\begin{theorem}}
\newcommand{\et}{\end{theorem}}
\newcommand{\bl}{\begin{lemma}}
\newcommand{\el}{\end{lemma}}
\newcommand{\bp}{\begin{proposition}}
\newcommand{\ep}{\end{proposition}}
\newcommand{\bcor}{\begin{corollary}}
\newcommand{\ecor}{\end{corollary}}
\newcommand{\br}{\begin{remark}\rm}
\newcommand{\er}{\end{remark}}
\newcommand{\bcon}{\begin{conjecture}}
\newcommand{\econ}{\end{conjecture}}
\newcommand{\bex}{\begin{exercise}}
\newcommand{\eex}{\end{exercise}}
\newcommand{\bcou}{\begin{counterex}}
\newcommand{\ecou}{\end{counterex}}

\newenvironment{Proof}[1][]{\noi\textbf{Proof #1}}{\QED}
\newcommand{\bpro}{\begin{Proof}}
\newcommand{\epro}{\end{Proof}}

\newcommand{\be}{\begin{equation}}
\newcommand{\ee}{\end{equation}}
\newcommand{\ba}{\begin{array}}
\newcommand{\ea}{\end{array}}
\newcommand{\bc}{\be\begin{array}{r@{\,}c@{\,}l}}
\newcommand{\ec}{\end{array}\ee}
\newcommand{\bac}{\begin{array}{r@{\,}c@{\,}l}}

\newcommand{\ga}{\gamma}

\newcommand{\de}{\delta}

\newcommand{\La}{\Lambda}
\newcommand{\sig}{\sigma}

\newcommand{\om}{\omega}
\newcommand{\Om}{\Omega}

\newcommand{\si}{\ensuremath{\sigma}}

\newcommand{\Ci}{{\cal C}}
\newcommand{\Di}{{\cal D}}

\newcommand{\Fi}{{\cal F}}
\newcommand{\Gi}{{\cal G}}
\newcommand{\Hi}{{\cal H}}

\newcommand{\Pc}{{\cal P}}

\newcommand{\Ti}{{\cal T}}

\newcommand{\E}{{\mathbb E}}
\newcommand{\N}{{\mathbb N}}
\renewcommand{\P}{{\mathbb P}}

\newcommand{\R}{{\mathbb R}}
\newcommand{\T}{{\mathbb T}}
\newcommand{\U}{{\mathbb U}}

\newcommand{\Z}{{\mathbb Z}}

\newcommand{\li}{\langle}
\newcommand{\re}{\rangle}

\newcommand{\desd}{\ensuremath{\Leftrightarrow}}
\newcommand{\volgt}{\ensuremath{\Rightarrow}}

\newcommand{\sub}{\subset}
\newcommand{\beh}{\backslash}

\newcommand{\asto}[1]{\underset{{#1}\to\infty}{\longrightarrow}}
\newcommand{\Asto}[1]{\underset{{#1}\to\infty}{\Longrightarrow}}

\newcommand{\ov}{\overline}
\newcommand{\un}{\underline}

\newcommand{\pa}{\partial}

\newcommand{\di}{\mathrm{d}}

\newcommand{\wurz}{\varnothing}

\newcommand{\ibf}{\mathbf{i}}
\newcommand{\jbf}{\mathbf{j}}
\newcommand{\kbf}{\mathbf{k}}

\newcommand{\ch}{\check}

\begin{document}

%numbering formulas within sections
\makeatletter\@addtoreset{equation}{section}
\makeatother\def\theequation{\thesection.\arabic{equation}} 

%alternative layout for enumerate lists.
\renewcommand{\labelenumi}{{\rm (\roman{enumi})}}

\title{A new characterization of endogeny}
\author{
Tibor~Mach
\footnote{The Czech Academy of Sciences,
Institute of Information Theory and Automation.
Pod vod\'arenskou v\v e\v z\' i 4,
18208 Praha 8,
Czech Republic;
swart@utia.cas.cz}
\and
Anja~Sturm
\footnote{Institute for Mathematical Stochastics,
Georg-August-Universit\"at G\"ottingen,
Goldschmidtstr.~7,
37077 G\"ottingen,
Germany;
asturm@math.uni-goettingen.de}
\and
Jan~M.~Swart$\;^\ast$
}

\date{\today}

\maketitle

\begin{abstract}\noi
Aldous and Bandyopadhyay have shown that each solution to a recursive
distributional equation (RDE) gives rise to a recursive tree process (RTP),
which is a sort of Markov chain in which time has a tree-like structure and in
which the state of each vertex is a random function of its descendants. If the
state at the root is measurable with respect to the sigma field generated by
the random functions attached to all vertices, then the RTP is said to be
endogenous. For RTPs defined by continuous maps, Aldous and Bandyopadhyay
showed that endogeny is equivalent to bivariate uniqueness, and they asked if
the continuity hypothesis can be removed. We introduce a higher-level RDE
that through its $n$-th moment measures contains all $n$-variate RDEs. We show
that this higher-level RDE has minimal and maximal fixed points with respect
to the convex order, and that these coincide if and only if the corresponding
RTP is endogenous. As a side result, this allows us to answer the question
of Aldous and Bandyopadhyay positively.
\end{abstract}
\vspace{.5cm}

\noi
{\it Acknowledgement.} Work sponsored by grant 16-15238S of the Czech Science
Foundation (GA CR). %Kotecky-Swart 16-18

%60K35   	Interacting random processes; statistical mechanics type
%               models; percolation theory
%
%60J05   	Discrete-time Markov processes on general state spaces
%82C22   	Interacting particle systems
%60J80   	Branching processes (Galton-Watson, birth-and-death, etc.)

{\setlength{\parskip}{-7pt}\tableofcontents}
\newpage

\section{Recursive distributional equations}

Let $S$ be a Polish space and for $k\geq 1$, let $S^k$ denote the space of all
ordered sequences $(x_1,\ldots,x_k)$ of elements of $S$. Let $\Pc(S)$ denote
the space of all probability measures on $S$, equipped with the topology of
weak convergence and its associated Borel-\si-field. A measurable map
$g:S^k\to S$ gives rise to a measurable map $\ch g:\Pc(S)^k\to\Pc(S)$ defined
as
%A proof that this map is indeed measurable can be found at the end of the
%paper, after the \end{document}.
\be\label{chg}
\ch g(\mu_1,\ldots,\mu_k):=\mu_1\otimes\cdots\otimes\mu_k\circ g^{-1},
\ee
where $\mu_1\otimes\cdots\otimes\mu_k$ denotes product measure and the
right-hand side of (\ref{chg}) is the image of this under the map $g$.
A more probabilistic way to express (\ref{chg}) is to say that if
$X_1,\ldots,X_k$ are independent random variables with laws
$\mu_1,\ldots,\mu_k$, then $g(X_1,\ldots,X_k)$ has law $\ch
g(\mu_1,\ldots,\mu_k)$. In particular, we let $T_g:\Pc(S)\to\Pc(S)$ denote the
map
\be
\label{T_gmap}
T_g(\mu):=\ch g(\mu,\ldots,\mu).
\ee
Note that $T_g$ is in general nonlinear, unless $k=1$.

With slight changes in the notation, the construction above works also for
$k=0$ and $k=\infty$. By definition, we let $S^0$ be a set containing a single
element, the empty sequence, and we let $S^\infty$ denote the space of all 
infinite sequences $(x_1,x_2,\ldots)$ of elements of $S$, equipped with the 
product topology and associated Borel-\si-field. It is well-known that if $S$ is 
Polish, then so are $S^k$ $(0\leq k\leq\infty)$ and $\Pc(S)$.

Write $\N:=\{0,1,2,\ldots\}$ and $\ov\N:=\N\cup\{\infty\}$. Let $\Gi$ be a
measurable space whose elements are measurable maps $g:S^k\to S$, where
$k=k_g\in\ov\N$ may depend on $g$, and let $\pi$ be a probability law on
$\Gi$.  Then under suitable technical assumptions (to be made precise in
  the next section)
\be\label{Tmap}
T(\mu):=\int_\Gi\pi(\di g)\,T_g(\mu)\qquad\big(\mu\in\Pc(S)\big)
\ee
defines a map $T:\Pc(S)\to\Pc(S)$. Equations of the form
\be\label{RDE}
T(\mu)=\mu
\ee
are called \emph{Recursive Distributional Equations} (RDEs). A nice collection
of examples of such RDEs arising in a variety of settings can be found in
\cite{AB05}. They include Galton-Watson branching processes and related random
trees, probabilistic analysis of algorithms as well as statistical physics
models.

\section{Recursive Tree Processes}

We now make our assumptions on the set $\Gi$ and probability measure
$\pi$ from (\ref{Tmap}) explicit. We fix a Polish space $\Om$, equipped with
the Borel \si-field, modeling some source of external randomness. We
let $\kappa:\Om\to\ov\N$ be measurable, set
$\Om_k:=\{\om\in\Om:\kappa(\om)=k\}$ $(k\in\ov\N)$, and let $\ga$ be a
function such that for each $k\in\ov\N$,
\be\label{gamma}
\Om_k\times S^k\ni(\om,x)\mapsto\ga[\om](x)\in S
\quad\mbox{is jointly measurable in $\om$ and $x$.}
\ee
Note that formally, $\ga$ is a
function from $\bigcup_{k\in\ov\N}\Om_k\times S^k$ into $S$.
We let $r$ be a probability measure on $\Om$.

The joint measurability of $\ga[\om](x)$ in $\om$ and $x$ implies in
particular that $\ga[\om]:S^{\kappa(\om)}\to S$ is measurable for each
$\om\in\Om$. We assume that the set $\Gi$ from (\ref{Tmap}) is given by
\be\label{Gidef}
\Gi:=\{\ga[\om]:\om\in\Om\},
\ee
which we equip with the final \si-field for the map $\om\mapsto\ga[\om]$. We
further assume that the probability measure $\pi$ from (\ref{Tmap}) is the image
of the probability measure $r$ on $\Om$ under this map.

Let $\li\mu,\phi\re:=\int_S\mu(\di x)\,\phi(x)$ denote the integral of a
bounded measurable function $\phi:S\to\R$ w.r.t.\ a measure $\mu$ on
$S$. Then, under the assumptions we have just made,
\bc\label{Tmudef2}
\dis\li T(\mu),\phi\re
&=&\dis\int_\Gi\!\pi(\di g)\,\li T_g(\mu),\phi\re
=\int_\Om\!r(\di\om)\,\li T_{\ga[\om]}(\mu),\phi\re\\[5pt]
&=&\dis\sum_{k\in\ov\N}\int_{\Om_k}\!r(\di\om)\,
\int_S\!\mu(\di x_1)\,\cdots\int_S\!\mu(\di x_k)\,
\phi\big(\ga[\om](x_1,\ldots,x_k)\big).
\ec
The joint measurability of $\ga[\om](x)$ in $\om\in\Om_k$ and $x\in S^k$
guarantees that $\li T_{\ga[\om]}(\mu),\phi\re$, which is defined by repeated
integrals over $S$, is measurable as a function of $\om$ and hence the
integral over $\Om$ is well-defined. Using this, one can check that our choice
of the \si-field on $\Gi$ guarantees that $\li T_g(\mu),\phi\re$ is measurable
as a function of $g$, and (\ref{Tmudef2}) defines a probability measure on $S$.

Starting from any countable set $\Gi$ of measurable maps $g:S^k\to S$ (where
$k=k_g$ may depend on $g$) and a probability law $\pi$ on $\Gi$, it is easy to
see that one can always construct $\Om$, $r$, and $\ga$ in terms of which
$\Gi$ and $\pi$ can then be constructed as above. For uncountable $\Gi$, the
construction above not only serves as a convenient technical set-up that
guarantees that the map in (\ref{Tmap}) is well-defined, but also has a
natural interpretation, with $\om$ playing the role of a source of external
randomness. This external randomness plays a natural role in Recursive
Tree Processes as introduced in \cite{AB05}, which we describe next.

For $d\in\N_+:=\{1,2,\ldots\}$, let $\T^d$ denote the space of all finite
words $\ibf=i_1\cdots i_t$ $(t\in\N)$ made up from the alphabet
$\{1,\ldots,d\}$, and define $\T^\infty$ similarly, using the alphabet
$\N_+$. Let $\wurz$ denote the word of length zero. We view $\T^d$ as a
tree with root $\wurz$, where each vertex $\ibf\in\T^d$ has $d$ children
$\ibf 1,\ibf 2,\ldots$, and each vertex $\ibf=i_1\cdots i_t$ except the root
has precisely one ancestor $i_1\cdots i_{t-1}$. If $\ibf, \jbf \in
\T^d$ with  $\ibf=i_1\cdots i_s$ and $\jbf=j_1\cdots j_t$, then we define the
concatenation $\ibf \jbf \in \T^d$ by $\ibf \jbf =i_1\cdots i_s j_1\cdots j_t$.
We denote the length of a word $\ibf=i_1\cdots i_t$ by $|\ibf|:=t$ and set
$\T^d_{(t)}:=\{\ibf\in\T^d:|\ibf|<t\}$. For any subtree $\U\sub\T$, we let
$\pa\U:=\{\ibf=i_1\cdots i_{|\ibf|} \in\T^d:
i_1\cdots i_{|\ibf|-1}\in\U,\ \ibf\not\in\U\}$ denote the outer
boundary of $\U$. In particular, $\pa\T^d_{(t)}=\{\ibf\in\T^d:|\ibf|=t\}$ is the
set of all vertices at distance $t$ from the root.

Fix some $d\in\ov\N_+:=\N_+\cup\{\infty\}$ such
that $\kappa(\om)\leq d$ for all $\om\in\Om,$ 
and to simplify notation write
$\T:=\T^d$. Let $(\om_\ibf)_{\ibf\in\T}$ be an i.i.d.\ collection
of $\Om$-valued random variables with common law $r$.
Fix $t\geq 1$ and $\mu\in\Pc(S)$, and let
$(X_\ibf)_{\ibf\in\pa \T_{(t)}}$ be a collection of $S$-valued random variables
such that
\be\label{finRTP1}
(X_\ibf)_{\ibf\in\pa \T_{(t)}}\mbox{ are i.i.d.\ with common law }\mu
\mbox{ and independent of }(\om_\ibf)_{\ibf\in\T_{(t)}}.
\ee
Define $(X_\ibf)_{\ibf\in\T_{(t)}}$ inductively by\footnote{Here and in similar
  formulas to come, it is understood that the notation should be suitably
  adapted if $\kappa(\om_\ibf)=\infty$, e.g., in this place,
  $X_\ibf=\ga [\om_\ibf]\big(X_{\ibf 1},X_{\ibf 2},\ldots\big)$.}
\be\label{finRTP2}
X_\ibf:=\ga[\om_\ibf](X_{\ibf 1},\ldots,X_{\ibf\kappa(\om_\ibf)})
\ee
Then it is easy to see that for each $1\leq s\leq t$,
\be
(X_\ibf)_{\ibf\in\pa\T_{(s)}}\mbox{ are i.i.d.\ with common law }T^{t-s}(\mu)
\mbox{ and independent of }(\om_\ibf)_{\ibf\in\T_{(s)}},
\ee
where $T^n$ denotes the $n$-th iterate of the map in (\ref{Tmap}).
Also, $X_\wurz$ (the state at the root) has law $T^t(\mu)$.
If $\mu$ is a solution of the RDE (\ref{RDE}), then, by
Kolmogorov's extension theorem  
there exists a collection
$(\om_\ibf,X_\ibf)_{\ibf\in\T}$ of random variables whose joint law is
uniquely characterized by the following requirements:
\be\ba{rl}\label{RTP}
{\rm(i)}&\mbox{$(\om_\ibf)_{\ibf\in\T}$ is an i.i.d.\ collection of
  $\Om$-valued r.v.'s with common law $r$,}\\[5pt]
{\rm(ii)}&\mbox{for each $t\geq 1$, the $(X_\ibf)_{\ibf\in{\pa\T_{(t)}}}$ are
  i.i.d.\ with common law $\mu$}\\
&\mbox{and independent of $(\om_\ibf)_{\ibf\in\T_{(t)}}$,}\\[5pt]
{\rm(iii)}&\dis X_\ibf:=\ga[\om_\ibf](X_{\ibf 1},\ldots,X_{\ibf\kappa(\om_\ibf)})
\qquad(\ibf\in\T).
\ec
We call such a collection $(\om_\ibf,X_\ibf)_{\ibf\in\T}$ a \emph{Recursive Tree
  Process} (RTP) corresponding to the map $\ga$ and the solution $\mu$ of the RDE (\ref{RDE}).
  We can think of an RTP as a generalization of a stationary and time reversed Markov
chain, where the time index set $\T$ has a tree structure and time flows in the
direction of the root. In each step, the new value $X_\ibf$ is a
function of the previous values $X_{\ibf1},X_{\ibf2}\ldots,$
plus some independent randomness, represented by the random variables 
$(\om_\ibf)_{\ibf\in\T}.$
Following \cite[Def~7]{AB05}, we say that the RTP
corresponding to a solution $\mu$ of the RDE (\ref{RDE}) is \emph{endogenous}
if $X_\wurz$ is measurable w.r.t.\ the \si-field generated by the random
variables $(\om_\ibf)_{\ibf\in \T}$.

Endogeny is somewhat similar to pathwise uniqueness of stochastic differential
equations, in the sense that it asks whether given $(\om_\ibf)_{\ibf\in\T}$,
there always exists a ``strong solution'' $(X_\ibf)_{\ibf\in\T}$ on the same
probability space, or whether on the other hand additional randomness is
needed to construct $(X_\ibf)_{\ibf\in\T}$. Since for each $t\geq 1$,
$X_\wurz$ is a function of $(\om_\ibf)_{\ibf\in\T_{(t)}}$ and the ``boundary
conditions'' $(X_\ibf)_{\ibf\in{\pa\T_{(t)}}}$, endogeny says that in a certain
almost sure sense, the effect of the boundary conditions disappears as
$t\to\infty$. Nevertheless, endogeny does not imply uniqueness of solutions to
the RDE (\ref{RDE}). Indeed, it is possible for a RDE to have several
solutions, while some of the corresponding RTPs are endogenous and others are
not. In the special case that $\T=\T^1$, an RTP is a time reversed stationary
Markov chain $\ldots,X_{11},X_1,X_\wurz$ generated by i.i.d.\  random
variables $\ldots,\om_{11},\om_{1}$. In this context, equivalent formulations of
endogeny have been investigated in the literature, see for example \cite{BL07}
who point back to \cite{Ros59}. Endogeny also plays a role, for example, in
the coupling from the past algorithm by Propp and Wilson \cite{PW96}. We point
to Lemma 14 and 15 of Section 2.6 in \cite{AB05} for an analogous statement on 
tree-structured coupling from the past.

Endogeny of RTPs is related to, but not the same as triviality of the
 tail-\si-field as defined for infinite volume Gibbs measures. Indeed,
if $(\om_\ibf,X_\ibf)_{\ibf\in\T}$ is an RTP, then the law of
$(X_\ibf)_{\ibf\in\T}$ is a Gibbs measure on $S^\T$. The \emph{tail-\si-field}
of such a Gibbs measure is defined as 
\be
\Ti:=\bigcap_{t\geq 1}\sig\big((X_\ibf)_{\ibf\in\T\beh\T_{(t)}}\big).
\ee
It is known that if $(\om_\ibf,X_\ibf)_{\ibf\in\T}$ is endogenous, then the
tail-\si-field of $(X_\ibf)_{\ibf\in\T}$ is trivial \cite[Prop.~1]{Ant06}, but
the converse implication does not hold \cite[Example~1]{Ant06}. It is known
that triviality of the tail-\si-field is equivalent to nonreconstructability
in information theory \cite[Prop.~15]{Mos01}, and also to extremality of
$(X_\ibf)_{\ibf\in\T}$ as a Gibbs measure \cite[Section~7.1]{Geo11}.

\section{The n-variate RDE}

Let $g:S^k\to S$ with $k\geq 0$ be a measurable map and let $n\geq 1$ be an
integer. We can naturally identify the space $(S^n)^k$ with the space of all
$n\times k$ matrices $x=(x_i^j)_{i=1,\ldots,k}^{j=1,\ldots,n}$. We let
$x^j:=(x^j_1,\ldots,x^j_k)$ and $x_i=(x^1_i,\ldots,x^n_i)$ denote the rows and
columns of such a matrix, respectively. With this notation, we define an
\emph{$n$-variate map} $g^{(n)}:(S^n)^k\to S^n$ by
\be\label{nvar}
g^{(n)}\big(x):=\big(g(x^1),\ldots,g(x^n)\big)
\qquad\big(x\in(S^n)^k\big).
\ee
This notation is easily generalized to $k=\infty$ or $n=\infty$, or both.
The map $g^{(n)}$ describes $n$ systems that are coupled in such a
way that the same map $g$ is applied to each system. We will be interested in
the \emph{$n$-variate map} (compare (\ref{Tmap}))
\be
\label{nvarTmap}
T^{(n)}(\nu):=\int_\Gi\pi(\di g)\,T_{g^{(n)}}(\nu)
\qquad\big(\nu\in\Pc(S^n)\big).
\ee
and the corresponding \emph{$n$-variate RDE} (compare (\ref{RDE}))
\be\label{nvarRDE}
T^{(n)}(\nu)=\nu.
\ee
If $\ga$ satisfies (\ref{gamma}), then the same is true for $\ga^{(n)}$ so
$T^{(n)}(\nu)$ is well-defined. The maps $T^{(n)}$ are consistent in the
following sense. Let $\nu|_{\{i_1,\ldots,i_m\}}$ denote the mar\-gi\-nal of
$\nu$ with respect to the coordinates $i_1,\ldots,i_m$, i.e., the image of
$\nu$ under the projection $(x_1,\ldots,x_n)\mapsto(x_{i_1},\ldots,x_{i_m})$. Then
\be
T^{(n)}(\nu)\big|_{\{i_1,\ldots,i_m\}}
=T^{(m)}\big(\nu\big|_{\{i_1,\ldots,i_m\}}\big).
\ee
In particular, if $\nu$ solves the $n$-variate RDE (\ref{nvarRDE}), then
its one-dimensional marginals $\nu|_{\{m\}}$ $(1\leq m\leq n)$ solve the RDE
(\ref{RDE}). For any $\mu\in\Pc(S)$, we let
\be
\Pc(S^n)_\mu
:=\big\{\nu\in\Pc(S^n):\nu|_{\{m\}}=\mu\ \forall 1\leq m\leq n\big\}
\ee
denote the set of probability measures on $S^n$ whose one-dimensional
marginals are all equal to $\mu$. We also let $\Pc_{\rm sym}(S^n)$ denote the
space of all probability measures on $S^n$ that are symmetric with respect to
permutations of the coordinates $\{1,\ldots,n\}$, and denote
$\Pc_{\rm sym}(S^n)_\mu:=\Pc_{\rm sym}(S^n)\cap\Pc(S^n)_\mu$.
It is easy to see that $T^{(n)}$ maps $\Pc_{\rm sym}(S^n)$ into itself.
If $\mu$ solves the RDE (\ref{RDE}), then $T^{(n)}$ also maps
$\Pc_{\rm sym}(S^n)_\mu$ into itself.

Given a measure $\mu\in\Pc(S)$, we define $\ov\mu^{(n)}\in\Pc(S^n)$
by
\be\label{diag}
\ov\mu^{(n)}:=\P\big[(X,\ldots,X)\in\,\cdot\,\big]
\quad\mbox{where $X$ has law }\mu.
\ee
We will prove the following theorem, which is similar to
\cite[Thm~11]{AB05}. The main improvement compared to the latter is that the
implication (ii)$\volgt$(i) is shown without the additional assumption that
$T^{(2)}$ is continuous with respect to weak convergence, solving Open
Problem~12 of \cite{AB05}. We have learned that this problem has been solved
before using an argument from \cite{BL07}, although its solution has not
been published. We refer to Appendix~\ref{A:problem} for a comparison of our
solution and this other solution. Below, $\Rightarrow$ denotes weak
convergence of probability measures.

\bt[Endogeny and bivariate uniqueness]
Let\label{T:bivar} $\mu$ be a solution to the RDE (\ref{RDE}). Then the
following statements are equivalent.
\begin{enumerate}
\item The RTP corresponding to $\mu$ is endogenous.
\item The measure $\ov\mu^{(2)}$ is the unique fixed point of $T^{(2)}$ in the
  space $\Pc_{\rm sym}(S^2)_\mu$.
\item $\dis(T^{(n)})^t(\nu)\Asto{t}\ov\mu^{(n)}$ for all $\nu\in\Pc(S^n)_\mu$
  and $n\in\ov\N_+$.
\end{enumerate}
\et

\section{The higher-level RDE}\label{S:hilev}

In this section we introduce a higher-level map $\ch T$ that through its
$n$-th moment measures contains all $n$-variate maps (Lemma~\ref{L:momeas}
below). In the next section, we will use this higher-level map to give a short
and elegant proof of Theorem~\ref{T:bivar}. We believe the methods of the
present section to be of wider interest. In particular, in future work we plan
to use them to study iterates of the $n$-variate maps for a non-endogenous RTP
related to systems with cooperative branching.

Let $\xi$ be a random probability measure on $S$, i.e., a $\Pc(S)$-valued
random variable, and let $\rho\in\Pc(\Pc(S))$ denote the law of $\xi$.
Conditional on $\xi$, let $X^1,\ldots,X^n$ be independent
with law $\xi$. Then (see Lemma~\ref{L:mommes} below)
\be\label{momdef}
\rho^{(n)}:=\P\big[(X^1,\ldots,X^n)\in\,\cdot\,\big]
=\E\big[\underbrace{\xi\otimes\cdots\otimes\xi}_{\mbox{$n$ times}}]
\ee
is called the \emph{$n$-th moment measure} of $\xi$. Here, the expectation of
a random measure $\xi$ on $S$ is defined in the usual way, i.e., $\E[\xi]$ is
the deterministic measure defined by
$\int\!\phi\,\di\E[\xi]:=\E[\int\!\phi\,\di\xi]$ for any bounded measurable
$\phi:S\to\R$. A similar definition applies to measures on $S^n$.
With slight changes in the notation, $\rho^{(n)}$ can also be defined for
$n=\infty$.

We observe that $\rho^{(n)}\in\Pc_{\rm sym}(S^n)$ for each
$\rho\in\Pc(\Pc(S))$ and $n\in\ov\N_+$. By De Finetti's theorem, for
$n=\infty$ the converse implication also holds. Indeed,
$\Pc_{\rm sym}(S^\infty)$ is the space of exchangeable probability measures on
$S^\infty$ and De Finetti says that each element of $\Pc_{\rm sym}(S^\infty)$
is of the form $\rho^{(\infty)}$ for some $\rho\in\Pc(\Pc(S))$.
Thus, we have a natural identification $\Pc_{\rm
  sym}(S^\infty)\cong\Pc(\Pc(S))$ and through this identification the map
$T^{(\infty)}:\Pc_{\rm sym}(S^\infty)\to\Pc_{\rm sym}(S^\infty)$ corresponds to
a map on $\Pc(\Pc(S))$. Our next aim is to identify this map.

Let $\ch T:\Pc(\Pc(S))\to\Pc(\Pc(S))$ be given by
\be
\label{higherTmap}
\ch T(\rho):=\int_\Gi\pi(\di g)\, T_{\ch g}(\rho)
\qquad\big(\rho\in\Pc(\Pc(S))\big),
\ee
where for any $g:S^k\to S$, the map $\ch g:\Pc(S)^k\to\Pc(S)$ is defined as in
(\ref{chg}). If $\ga$ satisfies (\ref{gamma}), then the same is true for
$\ch\ga$ so $\ch T(\rho)$ is well-defined.
%A proof of this claim can be found at the end of the
%paper, after the \end{document}.
Note that $T_{\ch g}(\rho)={\ch{\ch g}}(\rho,\ldots,\rho)$ by
(\ref{T_gmap}). We call $\ch T$ the \emph{higher-level map}, which gives rise
to the \emph{higher-level RDE}
\be\label{hlRDE}
\ch T(\rho)=\rho.
\ee
The following lemma shows that $\ch T$ is the map corresponding to
$T^{(\infty)}$ we were looking for. More generally, the lemma links $\ch T$ to
the $n$-variate maps $T^{(n)}$.

\bl[Moment measures]
Let\label{L:momeas} $n\in\ov\N_+$ and let $T^{(n)}$ and $\ch T$ be defined as in
(\ref{nvarTmap}) and (\ref{higherTmap}). Then the $n$-th moment measure of $\ch
T(\rho)$ is given by
\be\label{momeas}
\ch T(\rho)^{(n)}=T^{(n)}(\rho^{(n)})\qquad\big(\rho\in\Pc(\Pc(S))\big).
\ee
\el

Lemma~\ref{L:momeas} implies in particular that if $\rho$
solves the higher-level RDE (\ref{hlRDE}), then its first moment measure
$\rho^{(1)}$ solves the original RDE (\ref{RDE}). For any $\mu\in\Pc(S)$, we
let
\be\label{submu}
\Pc(\Pc(S))_\mu:=\big\{\rho\in\Pc(\Pc(S)):\rho^{(1)}=\mu\big\}
\ee
denote the set of all $\rho$ whose first moment measure is $\mu$. Note that
$\rho\in\Pc(\Pc(S))_\mu$ implies $\rho^{(n)}\in\Pc_{\rm sym}(S^n)_\mu$ for
each $n\geq 1$.

We equip $\Pc(\Pc(S))$ with the \emph{convex order}. By Theorem~\ref{T:Stras}
in Appendix~\ref{A:cv}, two measures $\rho_1,\rho_2\in\Pc(\Pc(S))$ are ordered
in the convex order, denoted $\rho_1\leq_{\rm cv}\rho_2$, if and only if there
exists an $S$-valued random variable $X$ defined on some probability space
$(\Om,\Fi,\P)$ and sub-\si-fields $\Fi_1\sub\Fi_2\sub\Fi$ such that
$\dis\rho_i=\P\big[\P[X\in\,\cdot\,|\Fi_i]\in\,\cdot\,\big]$ $(i=1,2)$.  It is
not hard to see that $\Pc(\Pc(S))_\mu$ has a minimal and maximal element
w.r.t.\ the convex order. For any $\mu\in\Pc(S)$, let us define
\be
\ov\mu:=\P\big[\de_X\in\,\cdot\,\big]
\quad\mbox{where $X$ has law }\mu.
\ee
Clearly $\de_\mu,\ov\mu\in\Pc(\Pc(S))_\mu$. Moreover (as will be proved in
Section~\ref{S:other} below)
\be\label{minmax}
\de_\mu\leq_{\rm cv}\rho\leq_{\rm cv}\ov\mu
\quad\mbox{for all}\quad
\rho\in\Pc(\Pc(S))_\mu.
\ee
In line with notation that has already been introduced in (\ref{diag}),
the $n$-th moment measures of $\de_\mu$ and $\ov\mu$ are given by
\be\label{ringmom}
\de^{(n)}_\mu=\P\big[(X_1,\ldots,X_n)\in\,\cdot\,\big]
\quand
\ov\mu^{(n)}=\P\big[(X,\ldots,X)\in\,\cdot\,\big],
\ee
where $X_1,\ldots,X_n$ are i.i.d.\ with common law $\mu$ and $X$ has law
$\mu$. The following proposition says that the higher-level RDE (\ref{hlRDE})
has a minimal and maximal solution with respect to the convex order.

\bp[Minimal and maximal solutions]
The\label{P:minmax} map $\ch T$ is monotone w.r.t.\ the convex order.
Let $\mu$ be a solution to the RDE (\ref{RDE}).
Then $\ch T$ maps $\Pc(\Pc(S))_\mu$ into itself.
There exists a unique $\un\mu\in\Pc(\Pc(S))_\mu$
such that
\be\label{ringdef}
\ch T^t(\de_\mu)\Asto{t}\un\mu,
\ee
where $\Rightarrow$ denotes weak convergence of measures on $\Pc(S)$, equipped
with the topology of weak convergence. The measures $\un\mu$ and $\ov\mu$
solve the higher-level RDE (\ref{hlRDE}), and any $\rho\in\Pc(\Pc(S))_\mu$
that solves the higher-level RDE (\ref{hlRDE}) must satisfy
\be\label{sandwich}
\un\mu\leq_{\rm cv}\rho\leq_{\rm cv}\ov\mu.
\ee
\ep

Since $\un\mu$ and $\ov\mu$ solve the higher-level RDE (\ref{hlRDE}),
there exist RTPs corresponding to $\un\mu$ and $\ov\mu$. The following
proposition gives an explicit description of these higher-level RTPs.

\bp[Higher-level RTPs]
Let\label{P:hlRTP} $(\om_\ibf,X_\ibf)_{\ibf\in\T}$ be an RTP corresponding to a
solution $\mu$ of the RDE (\ref{RDE}). Set
\be\label{xidef}
\xi_\ibf:=\P\big[X_\ibf\in\,\cdot\,|\,(\om_{\ibf\jbf})_{\jbf\in\T}\big].
\ee
Then $(\om_\ibf,\xi_\ibf)_{\ibf\in\T}$ is an RTP corresponding to the
map $\ch \ga$ and the solution $\un\mu$ of the higher-level RDE
(\ref{hlRDE}). Also, $(\om_\ibf,\de_{X_\ibf})_{\ibf\in\T}$ is an RTP
corresponding to  the map $\ch\ga$ and  $\ov\mu$.
\ep

In general, we can interpret the higher-level map $\ch T$ as follows. Fix
$t\geq 1$, and let $(X_\ibf,Y_\ibf)_{\ibf\in\pa \T_{(t)}}$ be
i.i.d.\ random variables, independent of $(\om_\ibf)_{\ibf\in
  \T_{(t)}}$, where the $X_\ibf$'s take values in $S$ and the
$Y_\ibf$'s take values in some arbitrary measurable space. Define
$(X_\ibf)_{\ibf\in\T_{(t)}}$ inductively as in (\ref{finRTP2}) and for
$\ibf\in\T_{(t)}\cup\pa\T_{(t)}$, let $\xi_\ibf$ denote the conditional law of
$X_\ibf$ given $(\om_{\ibf\jbf})_{\ibf\jbf\in\T_{(t)}}$ and
  $(Y_{\ibf\jbf})_{\ibf\jbf\in\pa\T_{(t)}}$. Then the
$(\xi_\ibf)_{\ibf\in\pa\T_{(t)}}$ are i.i.d.\ with some common law
$\rho\in\Pc(\Pc(S))$. Using Lemma~\ref{L:Tch} below, it is not hard to see
that for each $1\leq s\leq t$,
\be
(\xi_\ibf)_{\ibf\in\pa\T_{(s)}}\mbox{ are i.i.d.\ with common law }\ch T^{t-s}(\rho)
\mbox{ and independent of }(\om_\ibf)_{\ibf\in\T_{(s)}}.
\ee
Let $\mu$ denote the common law of the random variables
$(X_\ibf)_{\ibf\in\pa\T_{(t)}}$. We think of the random variables
$(Y_\ibf)_{\ibf\in\pa\T_{(t)}}$ as providing extra information about the
$(X_\ibf)_{\ibf\in\pa\T_{(t)}}$. The convex order measures how much extra
information we have. For $\xi_\ibf=\de_{X_\ibf}$, we have perfect information,
while on the other hand for $\xi_\ibf=\mu$ the $Y_\ibf$'s provided no extra
information. A solution to the higher-level RDE gives rise to a higher-level
RTP that can be interpreted as a normal (low-level) RTP
$(\om_\ibf,X_\ibf)_{\ibf\in\T}$ in which we have extra information about the
$(X_\ibf)_{\ibf\in\T}$. The solutions $\un\mu$ and $\ov\mu$ of the
higher-level RDE correspond to minimal and maximal knowledge about $X_\ibf$,
respectively, where either we know only $(\om_{\ibf\jbf})_{\jbf\in\T}$, or we
have full knowledge about $X_\ibf$.

We will derive Theorem~\ref{T:bivar} from the following theorem, which is our
main result.

\bt[The higher-level RDE]
Let\label{T:bivar2} $\mu$ be a solution to the RDE (\ref{RDE}). Then the
following statements are equivalent.
\begin{enumerate}\addtolength\itemsep{3pt}
\item The RTP corresponding to $\mu$ is endogenous.
\item $\un\mu=\ov\mu$.
\item $\dis\ch T^t(\rho)\Asto{t}\ov\mu$ for all $\rho\in\Pc(\Pc(S))_\mu$.
\end{enumerate}
\et

\section{Proof of the main theorem}

In this section, we use Lemma~\ref{L:momeas} and Propositions~\ref{P:minmax}
and \ref{P:hlRTP} to prove Theorems~\ref{T:bivar} and \ref{T:bivar2}. We need
one more lemma.

\bl[Convergence in probability]
Let\label{L:asco} $(\om_\ibf,X_\ibf)_{\ibf\in\T}$ be an endogenous RTP
corresponding to a solution $\mu$ of the RDE (\ref{RDE}), and let
$(Y_\ibf)_{\ibf\in\T}$ be an independent i.i.d.\ collection of $S$-valued
random variables with common law $\mu$. For each $t\geq 1$, set
$X^t_\ibf:=Y_\ibf$ $(\ibf\in\pa\T_{(t)})$, and define $(X^t_\ibf)_{\ibf\in\T_{(t)}}$
inductively by
\be\label{Xtdef}
X^t_\ibf=\ga[\om_\ibf](X^t_{\ibf 1},\ldots,X^t_{\ibf\kappa(\om_\ibf)})
\qquad(\ibf\in\T_{(t)}).
\ee
Then
\be\label{toP}
X^t_\wurz\asto{t}X_\wurz\quad\mbox{in probability.}
\ee
\el
\bpro
The argument is basically the same as in the proof of \cite[Thm~11~(c)]{AB05},
but for completeness, we give it here. Let $f,g:S\to\R$ be bounded and
continuous and let $\Fi_t$ resp.\ $\Fi_\infty$ be the \si-fields generated by
$(\om_\ibf)_{\ibf\in \T_{(t)}}$ resp.\ $(\om_\ibf)_{\ibf\in \T}$. Since
$X_\wurz$ and $X^t_\wurz$ are conditionally independent and
identically distributed given $\Fi_t$,
\be\ba{l}
\dis\E\big[f(X_\wurz)g(X^t_\wurz)\big]
=\E\big[\E\big[f(X_\wurz)\,\big|\,\Fi_t\big]
\E\big[g(X^t_\wurz)\,\big|\,\Fi_t\big]\big]\\[5pt]
\dis\quad=\E\big[\E\big[f(X_\wurz)\,\big|\,\Fi_t\big]
\E\big[g(X_\wurz)\,\big|\,\Fi_t\big]\big]\\[5pt]
\dis\quad\asto{t}\E\big[\E\big[f(X_\wurz)\,\big|\,\Fi_\infty\big]
\E\big[g(X_\wurz)\,\big|\,\Fi_\infty\big]\big]
=\E\big[f(X_\wurz)g(X_\wurz)\big],
\ec
where we have used martingale convergence
and in the last step also endogeny.
Since this holds for arbitrary $f,g$, we conclude that the law of
$(X_\wurz,X^t_\wurz)$ converges weakly to the law of
$(X_\wurz,X_\wurz)$, which implies (\ref{toP}).
\QED
\epro

\bpro[of Theorem~\ref{T:bivar2}]
If the RTP corresponding to $\mu$ is endogenous, then the random variable
$\xi_\wurz$ defined in (\ref{xidef}) satisfies
$\xi_\wurz=\de_{X_\wurz}$. By Proposition~\ref{P:hlRTP},
$\xi_\wurz$ and $\de_{X_\wurz}$ have laws $\un\mu$ and $\ov\mu$,
respectively, so (i)$\volgt$(ii). Conversely, if (i) does not hold, then
$\xi_\wurz$ is with positive probability not a delta measure, so
(i)$\desd$(ii).

The implication (iii)$\volgt$(ii) is immediate from the definition of $\un\mu$
in (\ref{ringdef}). To get the converse implication,
we observe that by Proposition~\ref{P:minmax}, $\ch T$ is monotone with
respect to the convex order, so (\ref{minmax}) implies
\be\label{lowup}
\ch T^t(\de_\mu)\leq_{\rm cv}\ch T^t(\rho)\leq_{\rm cv}\ch T^t(\ov\mu)
\qquad(t\geq 0).
\ee
By Proposition~\ref{P:minmax}, $\ch T$ maps $\Pc(\Pc(S))_\mu$ into
  itself, so $\ch T^t(\rho)\in\Pc(\Pc(S))_\mu$ for each $t\geq 0$,
and hence by Lemma~\ref{L:Pmu} in Appendix~\ref{A:cv}, the measures $(\ch
T^t(\rho))_{t\geq 1}$ are tight. By Proposition~\ref{P:minmax}, the left-hand
side of (\ref{lowup}) converges weakly to $\un\mu$ as $t\to\infty$ while the
right-hand side equals $\ov\mu$ for each $t$, so we obtain that any
subsequential limit $\ch T^{t_n}(\rho)\Rightarrow\rho_\ast$ satisfies
$\un\mu\leq_{\rm cv}\rho_\ast\leq_{\rm cv}\ov\mu$. In particular, this shows
that (ii)$\volgt$(iii).
\QED \epro

\bpro[of Theorem~\ref{T:bivar}]
The implication (iii)$\volgt$(ii) is trivial. By Lemma~\ref{L:momeas} and the
fact that $\un\mu$ and $\ov\mu$ solve the higher-level RDE, we see that (ii)
implies $\un\mu^{(2)}=\ov\mu^{(2)}$. By Proposition~\ref{P:minmax},
$\un\mu\leq_{\rm cv}\ov\mu$. Now Lemma~\ref{L:rho12} from Appendix~\ref{A:cv}
shows that $\un\mu^{(2)}=\ov\mu^{(2)}$ and $\un\mu\leq_{\rm cv}\ov\mu$ imply
$\un\mu=\ov\mu$, so applying Theorem~\ref{T:bivar2} we obtain that
(ii)$\volgt$(i).

To complete the proof, we will show that (i)$\volgt$(iii). Let
$(\om_\ibf,X_\ibf)_{\ibf\in\T}$ be an RTP corresponding to $\mu$ and let
$(Y^1_\ibf,\ldots,Y^n_\ibf)_{\ibf\in\T}$ be an independent i.i.d.\ collection
of $S^n$-valued random variables with common law $\nu$. For each $t\geq 1$ and
$1\leq m\leq n$, set $X^{m,t}_\ibf:=Y^m_\ibf\quad(\ibf\in\pa\T_{(t)})$, and define
$(X^{m,t}_\ibf)_{\ibf\in\T_{(t)}}$ inductively as in (\ref{Xtdef}). Then
$(X^{1,t}_\wurz,\ldots X^{n,t}_\wurz)$ has law $(T^{(n)})^t(\nu)$,
and using endogeny, Lemma~\ref{L:asco} tells us that
\be
(X^{1,t}_\wurz,\ldots X^{n,t}_\wurz)
\asto{t}(X_\wurz,\ldots,X_\wurz)\quad\mbox{in probability.}
\ee
Since the right-hand side has law $\ov\mu^{(n)}$ (recall
  \eqref{ringmom}), this completes the proof. With a slight change of
notation, this argument also works for $n=\infty$.
\QED \epro

\section{Other proofs}\label{S:other}

In this section, we provide the proofs of Lemma~\ref{L:momeas} and
Propositions \ref{P:minmax} and \ref{P:hlRTP}, as well as formula
(\ref{minmax}). We start with some preliminary observations.

\bl[Moment measures]
Let\label{L:mommes} $X^1,\ldots,X^n$ be
$S$-valued random variables such that conditionally on some \si-field $\Hi$,
the $X^1,\ldots,X^n$ are i.i.d.\ with (random) law
$\P[X^j\in\,\cdot\,|\Hi]=\xi$ a.s.\ $(j=1,\ldots,n)$.
Let $\rho\in\Pc(\Pc(S))$ denote the law of $\xi$, i.e.,
\be
\rho=\P\big[\P[X^j\in\,\cdot\,|\Hi]\in\,\cdot\,\big]\qquad(j=1,\ldots,n).
\ee
Then
\be\label{mommes}
\rho^{(n)}=\P\big[(X^1,\ldots,X^n)\in\,\cdot\,\big]
=\E\big[\underbrace{\xi\otimes\cdots\otimes\xi}_{\mbox{$n$ times}}].
\ee
\el
\bpro
This follows by writing
\be\ba{l}
\dis\E\big[\prod_{i=1}^nf_i(X^i)\big]
=\E\big[\prod_{i=1}^n\E[f_i(X^i)\,|\,\Hi]\big]\\[5pt]
\dis\quad=\int\rho(\di\xi)\int_{S^n}\xi(\di x_1)\cdots\xi(\di x_n)\,
f_1(x)\cdots f_n(x_n).
\ec
for arbitrary bounded measurable $f_i:S\to\R$.
%Note that this also proves the second equality in (\ref{momdef}) if we let
%$\Hi$ be the \si-field generated by $\xi$.
\QED \epro

\bl[Higher-level map]
Let\label{L:Tch} $X_1,\ldots,X_k$ and $\Hi_1,\ldots,\Hi_k$ be $S$-valued
random variables and \si-fields, respectively, such that
$(X_1,\Hi_1),\ldots,(X_k,\Hi_k)$ are i.i.d.\ Let
\be
\rho=\P\big[\P[X_i\in\,\cdot\,|\Hi_i]\in\,\cdot\,\big]\qquad(i=1,\ldots,k).
\ee
Let $\Hi_1\vee\cdots\vee\Hi_k$ denote the \si-field generated by
$\Hi_1,\ldots,\Hi_k$. Then, for each measurable $g:S^k\to S$,
\be\label{Tchg}
T_{ \ch g}(\rho)
=\P\big[\P[g(X_1,\ldots,X_k)\in\,\cdot\,|\Hi_1\vee\cdots\vee\Hi_k]
\in\,\cdot\,\big].
\ee
Similarly, if $\om$ is an independent $\Om$-valued random variable
with law $r$ and $\Fi$ is the \si-algebra generated by $\om$, then
\be\label{Tch}
\ch T(\rho)
=\P\big[\P[\ga[\om](X_1,\ldots,X_k)\in\,\cdot\,|
\Hi_1\vee\cdots\vee\Hi_k \vee\Fi]\in\,\cdot\,\big].
\ee
\el
\bpro
Let $\xi_i:=\P[X_i\in\,\cdot\,|\Hi_i]$. Then (i) $\xi_1,\ldots,\xi_k$ are
i.i.d.\ with common law $\rho$, and (ii) conditional on
$\Hi_1\vee\cdots\vee\Hi_k$, the random variables $X_1,\ldots,X_k$ are
independent with laws $\xi_1,\ldots,\xi_k$, respectively.
Now (ii) implies
\be
\P[g(X_1,\ldots,X_k)\in\,\cdot\,|\Hi_1\vee\cdots\vee\Hi_k]
=\ch g(\xi_1,\ldots,\xi_k)\quad{\rm a.s.}
\ee
and in view of (i), (\ref{Tchg}) follows.
Since $\om$ is independent of $(X_1,\Hi_1),\ldots,(X_k,\Hi_k)$,
conditional on $\Hi_1\vee\cdots\vee\Hi_k\vee\Fi$, the random variables
$X_1,\ldots,X_k$ are again independent with laws $\xi_1,\ldots,\xi_k$,
respectively, and hence
\be
\P[\ga[\om](X_1,\ldots,X_k)\in\,\cdot\,|\Hi_1\vee\cdots\vee\Hi_k\vee\Fi]
=\ch\ga[\om](\xi_1,\ldots,\xi_k)\quad{\rm a.s.},
\ee
which implies (\ref{Tch}). With a slight change in
notation, these formulas hold also for $k=\infty$.
\QED \epro

\bpro[of formula (\ref{minmax})]
Let $\xi$ be a $\Pc(S)$-valued random variable with law $\rho$ and conditional
on $\xi$, let $X$ be an $S$-valued random variable with law $\xi$. Let $\Fi_0$
be the trivial \si-field, let $\Fi_1$ be the \si-field generated by $\xi$, and
let $\Fi_2$ be the \si-field generated by $\xi$ and $X$. Then
$\Fi_0\sub\Fi_1\sub\Fi_2$. Since $\rho^{(1)}=\mu$, the random variable $X$ has
law $\mu$. Now
\be\ba{l}
\P\big[\P[X\in\,\cdot\,|\,\Fi_0]\in\,\cdot\,\big]
=\P\big[\mu\in\,\cdot\,\big]=\de_\mu,\\[5pt]
\P\big[\P[X\in\,\cdot\,|\,\Fi_1]\in\,\cdot\,\big]
=\P\big[\xi\in\,\cdot\,\big]=\rho,\\[5pt]
\P\big[\P[X\in\,\cdot\,|\,\Fi_2]\in\,\cdot\,\big]
=\P\big[\de_X\in\,\cdot\,\big]=\ov\mu.
\ec
This proves that $\de_\mu\leq_{\rm cv}\rho\leq_{\rm cv}\ov\mu$.
\QED \epro

\bpro[of Lemma~\ref{L:momeas}]
Let $\xi_1,\ldots,\xi_k$ be i.i.d.\ with common law $\rho$ and conditional on 
$\xi_1,\ldots,\xi_k$, let $(X^j_i)^{j=1,\ldots,n}_{i=1,\ldots,k}$ be
independent $S$-valued random variables such that $X^j_i$ has
law $\xi_i$. Let $\Hi_i$ denote the \si-field generated by $\xi_i$.
Then $\rho=\P[\P[X^j_i\in\,\cdot\,|\Hi_i]\in\,\cdot\,]$ for each $i,j$.
By (\ref{mommes}),
\be
\rho^{(n)}=\P\big[(X^1_i,\ldots,X^n_i)\in\,\cdot\,\big]\qquad(i=1,\ldots,k).
\ee
Set $X_i:=(X^1_i,\ldots,X^n_i)$ and $X^j:=(X^j_1,\ldots,X^j_k)$. 
Since $X_1,\ldots,X_k$ are independent with law $\rho^{(n)}$,
\be\label{Tnrn}
T_{g^{(n)}}(\rho^{(n)})
=\P\big[g^{(n)}(X_1,\ldots,X_k)\in\,\cdot\,\big]
=\P\big[\big(g(X^1),\ldots,g(X^n)\big)\in\,\cdot\,\big].
\ee
Let $\Hi:=\Hi_1\vee\cdots\vee\Hi_k$. Since
$(X^j_1,\Hi_1),\ldots,(X^j_k,\Hi_k)$ are i.i.d.\ for each $j$, formula
(\ref{Tchg}) tells us that $T_{\ch g}(\rho)
=\P\big[\P[g(X^j)\in\,\cdot\,|\Hi]\in\,\cdot\,\big]$ $(j=1,\ldots,n)$.
Since conditionally on $\Hi$, the $g(X^1),\ldots,g(X^n)$ are i.i.d., formula
(\ref{mommes}) tells us that
\be
T_{\ch g}(\rho)^{(n)}=\P\big[\big(g(X^1),\ldots,g(X^n)\big)\in\,\cdot\,\big].
\ee
Combining this with (\ref{Tnrn}), we see that $T_{\ch
  g}(\rho)^{(n)}=T_{g^{(n)}}(\rho^{(n)})$ for each $g\in\Gi$. 
Now (\ref{momeas}) follows by integrating w.r.t.\ $\pi$.
\QED \epro

\bpro[of Propositions~\ref{P:minmax} and \ref{P:hlRTP}]
We first show that $\ch T$ is monotone w.r.t.\ the convex order. Let
$\rho_1,\rho_2\in\Pc(\Pc(S))_\mu$ satisfy $\rho_1\leq_{\rm cv}\rho_2$.
Then we can construct $S$-valued random variables
$X_1,\ldots,X_k$ as well as \si-fields $(\Hi^j_i)^{j=1,2}_{i=1,\ldots,k}$, such that
$(X_1,\Hi^1_1,\Hi^2_1),$ $\ldots,(X_k,\Hi^1_k,\Hi^2_k)$ are i.i.d.,
\be
\rho_j=\P\big[\P[X_i\in\,\cdot\,|\Hi^j_i]\in\,\cdot\,\big]
\qquad(i=1,\ldots,k,\ j=1,2),
\ee
and $\Hi^1_i\sub\Hi^2_i$ for all $i=1,\ldots,k$. Let $\om$ be an independent
$\Om$-valued random variable with law $r$ and let $\Fi$ be the
  \si-field generated by $\om$. Then (\ref{Tch}) says that
\be
\ch T(\rho_j)
=\P\big[\P[\ga [\om](X_1,\ldots,X_k)\in\,\cdot\,
|\Hi^j_1\vee\cdots\vee\Hi^j_k \vee\Fi]
\in\,\cdot\,\big]\qquad(j=1,2).
\ee
Since $\Hi^1_1\vee\cdots\vee\Hi^1_k \vee\Fi \sub\Hi^2_1\vee\cdots\vee\Hi^2_k \vee\Fi$, this
proves that $\ch T(\rho_1)\leq_{\rm cv}\ch T(\rho_2)$.

Let $\mu$ be a solution to the RDE (\ref{RDE}). Then by Lemma~\ref{L:momeas}
 for $\rho\in\Pc(\Pc(S))_\mu$ we have $\ch
T(\rho)^{(1)}=T(\rho^{(1)})=T(\mu)=\mu$, proving that $\ch T$ maps
$\Pc(\Pc(S))_\mu$ into itself. It will be convenient to combine the proof of
the remaining statements of Proposition~\ref{P:minmax} with the proof of
Proposition~\ref{P:hlRTP}. To check (as claimed in Proposition~\ref{P:hlRTP})
that $(\om_\ibf,\xi_\ibf)_{\ibf\in\T}$ is an RTP corresponding to the map
$\ch \gamma$ and  to $\un\mu$, we
need to check that:
\begin{enumerate}
\item The $(\om_\ibf)_{\ibf\in\T}$ are i.i.d.
\item For each $t\geq 1$, the $(\xi_\ibf)_{\ibf\in\pa\T_{(t)}}$ are i.i.d.\ with
  common law $\un\mu$\\ and independent of $(\om_\ibf)_{\ibf\in\T_{(t)}}$.
\item $\xi_\ibf= \ch \ga[\om_\ibf](\xi_{\ibf 1},\ldots, \xi_{\ibf \kappa(\om_\ibf)})$
  $(\ibf\in\T)$.
\end{enumerate}
Here (i) is immediate. 
Since $\xi_\ibf$ depends only on $(\om_{\ibf\jbf})_{\jbf\in\T}$, it
is also clear that the $(\xi_\ibf)_{\ibf\in\pa\T_{(t)}}$ are i.i.d.\ and
independent of $(\om_\ibf)_{\ibf\in\T_{(t)}}$. To see that their common law is
$\un\mu$, we may equivalently show that $\xi_\wurz$ has law $\un\mu$.
Thus, we are left with the task to prove (iii) and
\begin{enumerate}\addtocounter{enumi}{3}
\item $\P[\xi_\wurz\in\,\cdot\,]=\un\mu$.
\end{enumerate}
Let $\Fi^\ibf$ denote the \si-field generated by
$(\om_{\ibf\jbf})_{\jbf\in\T}$. Then, for any $\ibf\in\T$,
\be
\xi_\ibf=\P[X_\ibf\in\,\cdot\,|\Fi^\ibf]
=\P\big[\ga[\om_\ibf](X_{\ibf 1},\ldots, X_{\ibf \kappa(\om_\ibf)})
\in\,\cdot\,\big|\Fi^\ibf\big].
\ee
Conditional on $\Fi^\ibf$, the random variables $X_{\ibf 1},\ldots,X_{\ibf
 k_\ibf}$ are independent with respective laws $\xi_{\ibf 1},\ldots,\xi_{\ibf k_\ibf}$,
and hence $\ga[\om_\ibf](X_{\ibf 1},\ldots, X_{\ibf \kappa(\om_\ibf)})$
has law
$\ch \ga[\om_\ibf](\xi_{\ibf 1},\ldots, \xi_{\ibf \kappa(\om_\ibf)})$,
proving (iii).

To prove also (iv), we first need to prove (\ref{ringdef}) from 
Proposition~\ref{P:minmax}. Fix $t\geq 1$ and for $\ibf\in\T_{(t)}\cup\pa\T_{(t)}$, let
$\Fi^\ibf_t$ denote the \si-field generated by
$\{\om_{\ibf\jbf}:\jbf\in\T,\ |\ibf\jbf|<t\}$. In particular, if
$\ibf\in\pa\T_{(t)}$, then $\Fi^\ibf_t$ is the trivial \si-field. Set
\be
\xi^t_\ibf:=\P[X_\ibf\in\,\cdot\,|\Fi^\ibf_t]\qquad(\ibf\in\T_{(t)}\cup\pa\T_{(t)}).
\ee
In particular, $\xi^t_\ibf=\mu$ a.s.\ for $\ibf\in\pa\T_{(t)}$. Arguing as before,
we see that
\be
\xi^t_\ibf=\ch\ga[\om_\ibf]\big(\xi^t_{\ibf 1},\ldots,\xi^t_{\ibf \kappa(\om_\ibf)}\big)
\qquad(\ibf\in\T_{(t)}),
\ee
and hence
\be\label{Ttconv}
\ch T^t(\de_\mu)=\P[\xi^t_\wurz\in\,\cdot\,].
\ee
By martingale convergence,
\be\label{Tmart}
\xi^t_\wurz=\P[X_\wurz\in\,\cdot\,|\Fi^\wurz_t]
\asto{t}\P[X_\wurz\in\,\cdot\,|(\om_\ibf)_{\ibf\in\T}]=\xi_\wurz
\quad{\rm a.s.}
\ee
Combining this with (\ref{Ttconv}), we obtain (\ref{ringdef}) where $\un\mu$
is in fact the law of $\xi_\wurz$, proving (iv) as well.
This completes the proof that $(\om_\ibf,\xi_\ibf)_{\ibf\in\T}$ 
 is an RTP corresponding to  the map $\ch \gamma$ and  $\un\mu$.

The proof that  $(\om_\ibf,\de_{X_\ibf})_{\ibf\in\T}$ 
 is an RTP corresponding to  the map $\ch \gamma$ and  $\ov\mu$ is simpler. It is clear that (i) the
$(\om_\ibf)_{\ibf\in\T}$   
are i.i.d., and (ii) for each $t\geq 1$, the
$(\de_{X_\ibf})_{\ibf\in\pa\T_{(t)}}$ are i.i.d.\ with common law $\ov\mu$ and
independent of $(\om_\ibf)_{\ibf\in\T_{(t)}}$.  
To prove that also (iii)
$\dis\de_{X_\ibf}=\ch\ga[\om_\ibf]\big(\de_{X_{\ibf 1}},\ldots,\de_{X_{\ibf \kappa(\om_\ibf)}}\big)$
$(\ibf\in\T)$, it suffices to show that for any measurable $g:S^k\to S$,
\be\label{crwn}
\ch g(\de_{x_1},\ldots,\de_{x_k})=\de_{g(x_1,\ldots,x_k)}.
\ee
By definition, the left-hand side of this equation is the law of
$g(X_1,\ldots,X_k)$, where $X_1,\ldots,X_k$ are independent with laws
$\de_{x_1},\ldots,\de_{x_k}$, so the statement is obvious.

This completes the proof of Proposition~\ref{P:hlRTP}. Moreover, since the
marginal law of an RTP solves the corresponding RDE, our proof also shows that
the measures $\un\mu$ and $\ov\mu$ solve the higher-level RDE
(\ref{hlRDE}).

In view of this, to complete the proof of Proposition~\ref{P:minmax}, it
suffices to prove (\ref{sandwich}). If $\rho$ solves the higher-level RDE
(\ref{hlRDE}), then applying $\ch T^t$ to (\ref{minmax}), using the
monotonicity of $\ch T$ with respect to the convex order, we see that
$\ch T^t(\de_\mu)\leq_{\rm cv}\rho\leq_{\rm cv}\ov\mu$ for all $t$.
Letting $t\to\infty$, (\ref{sandwich}) follows.
\QED \epro

\appendix

\normalsize

\section{The convex order}\label{A:cv}

By definition, a \emph{$G_\de$-set} is a set that is a countable intersection
of open sets. By \cite[\S 6 No.~1, Theorem.~1]{Bou58}, for a metrizable space
$S$, the following statements are equivalent.
\begin{enumerate}
\item $S$ is Polish.
\item There exists a metrizable compactification $\ov S$ of $S$ such that $S$
 is a $G_\de$-subset of $\ov S$.
\item For each metrizable compactification $\ov S$ of $S$, $S$
 is a $G_\de$-subset of $\ov S$.
\end{enumerate}
Moreover, a subset $S'\sub S$ of a Polish space $S$ is Polish in the induced
topology if and only if $S'$ is a $G_\de$-subset of $S$.

Let $S$ be a Polish space. Recall that $\Pc(S)$ denotes the space of
probability measures on $S$, equipped with the topology of weak convergence.
In what follows, we fix a metrizable compactification $\ov S$ of $S$.
Then we can identify the space $\Pc(S)$ (including its topology) with the space
of probability measures $\mu$ on $\ov S$ such that $\mu(S)=1$. By Prohorov's
theorem, $\Pc(\ov S)$ is compact, so $\Pc(\ov S)$ is a metrizable
compactification of $\Pc(S)$. Recall the definition of $\Pc(\Pc(S))_\mu$ from 
(\ref{submu}).

\bl[Measures with given mean]
For\label{L:Pmu} any $\mu\in\Pc(S)$, the space $\Pc(\Pc(S))_\mu$ is compact.
\el
\bpro
Since any $\rho\in\Pc(\Pc(\ov S))$ whose first moment measure is $\mu$ must be
concentrated on $\Pc(S)$, we can identify $\Pc(\Pc(S))_\mu$ with the space of
probability measures on $\Pc(\ov S)$ whose first moment measure is $\mu$.
From this we see that $\Pc(\Pc(S))_\mu$ is a closed subset of $\Pc(\Pc(\ov
S))$ and hence compact.
\QED \epro

%It is also possible to prove Lemma~\ref{L:Pmu} directly, without appealing to
%compactifications. It is clear that $\Pc(\Pc(S))_\mu$ is closed, so it
%suffices to prove tightness. For every $\eps>0$ and $n\geq 0$ we can find a
%compact $K_{\eps,n}\sub S$ such that $\mu(S\beh K_{\eps,n})\leq\eps^24^{-n}$
%and hence, for any $\rho\in\Pc(\Pc(S))_\mu$,
%\[
%\rho\big(\{\mu\in\Pc(S):\mu(S\beh K_{\eps,n})>\eps 2^{-n}\}\big)\leq\eps 2^{-n}
%\]
%and as a result
%\[
%\rho\big(\{\mu\in\Pc(S):\mu(S\beh K_{\eps,n})>\eps 2^{-n}
%\mbox{ for some }n\geq 0\}\big)\leq\eps.
%\]
%Since
%\[
%C_\eps:=\big\{\mu\in\Pc(S):\mu(S\beh K_{\eps,n})\leq\eps 2^{-n}
%\ \forall n\geq 0\big\}
%\]
%is a compact subset of $\Pc(S)$, the tightness of $\Pc(\Pc(S))_\mu$ follows.

We let $\Ci(\ov S)$ denote the space of all continuous real functions on $\ov
S$, equipped with the supremumnorm, and we let $B(\ov S)$ denote the space of
bounded measurable real functions on $\ov S$. The following fact is well-known
(see, e.g., \cite[Cor~12.11]{Car00}).

\bl[Space of continuous functions]
$\Ci(\ov S)$\label{L:separ} is a separable Banach space.
\el

For each $f\in\Ci(\ov S)$, we define an affine function $l_f\in\Ci(\Pc(\ov
S))$ by $l_f(\mu):=\int f\,\di\mu$. The following lemma says that all
continuous affine functions on $\Pc(\ov S)$ are of this form.

\bl[Continuous affine functions]
A\label{L:conaf} function $\phi\in\Ci(\Pc(\ov S))$ is affine if and only if
$\phi=l_f$ for some $f\in\Ci(\ov S)$.
\el
\bpro
Let $\phi:\Pc(\ov S)\to\R$ be affine and continuous. Since $\phi$ is continuous,
setting $f(x):=\phi(\de_x)$ $(x\in\ov S)$ defines a continuous function $f:\ov
S\to\R$. Since $\phi$ is affine, $\phi(\mu)=l_f(\mu)$ whenever $\mu$ is a
finite convex combination of delta measures. Since such measures are dense in
$\Pc(\ov S)$ and $\phi$ is continuous, we conclude that $\phi=l_f$.
\QED \epro

\bl[Lower semi-continuous convex functions]
Let\label{L:convrep} $C\sub\Ci(\ov S)$ be convex, closed, and nonempty. Then
\be\label{convrep}
\phi:=\sup_{f\in C}\,l_f
\ee
defines a lower semi-continuous convex function
$\phi:\Pc(\ov S)\to(-\infty,\infty]$. Conversely, each such $\phi$ is of the
form (\ref{convrep}).
\el
\bpro
It is straightforward to check that (\ref{convrep}) defines a
lower semi-continuous convex function $\phi:\Pc(\ov S)\to(-\infty,\infty]$.
To prove that every such function is of the form (\ref{convrep}),
let $\Ci(\ov S)'$ denote the dual of the Banach space $\Ci(\ov S)$, i.e.,
$\Ci(\ov S)'$ is the space of all continuous linear forms $l:\Ci(\ov S)\to\R$.
We equip $\Ci(\ov S)'$ with the weak-$\ast$ topology, i.e., the weakest
topology that makes the maps $l\mapsto l(f)$ continuous for all $f\in\Ci(\ov
S)$. Then $\Ci(\ov S)'$ is a locally convex topological vector space
%defined by the system of seminorms $\{p_f\}$ defined as $p_f(l):=|l(f)|$
and by the Riesz-Markov-Kakutani representation theorem, we can view $\Pc(\ov
S)$ as a convex compact metrizable subset of $\Ci(\ov S)'$. Now any lower
semi-continuous convex function $\phi:\Pc(\ov S)\to(-\infty,\infty]$
can be extended to $\Ci(\ov S)'$ by putting $\phi:=\infty$ on the complement of
$\Pc(\ov S)$. Applying \cite[Thm~I.3]{CV77} we obtain that $\phi$ is the
supremum of all continuous affine functions that lie below it. By
Lemma~\ref{L:conaf}, we can restrict ourselves to continuous affine functions
of the form $l_f$ with $f\in\Ci(\ov S)$. It is easy to see that $\{f\in\Ci(\ov
S):l_f\leq\phi\}$ is closed and convex, proving that every lower
semi-continuous convex function $\phi:\Pc(\ov S)\to(-\infty,\infty]$ is of the
  form (\ref{convrep}).
\QED \epro

We define
\be
\Ci_{\rm cv}\big(\Pc(\ov S)\big)
:=\big\{\phi\in\Ci(\Pc(\ov S)):\phi\mbox{ is convex}\big\}
\ee
If two probability measures $\rho_1,\rho_2\in\Pc(\Pc(S))$ satisfy the
equivalent conditions of the following theorem, then we say that they are
ordered in the \emph{convex order}, and we denote this as $\rho_1\leq_{\rm
  cv}\rho_2$. The fact that $\leq_{\rm cv}$ defines a partial order will be
proved in Lemma~\ref{L:disdet} below. The convex order can be defined more
generally for $\rho_1,\rho_2\in\Pc(C)$ where $C$ is a convex space, but in the
present paper we will only need the case $C=\Pc(\ov S)$.

\bt[The convex order for laws of random probability measures]
Let\label{T:Stras} $S$ be a Polish space and let $\ov S$ be a metrizable
compactification of $S$. Then, for $\rho_1,\rho_2\in\Pc(\Pc(S))$, the following
statements are equivalent.
\begin{enumerate}
\item $\dis\int\phi\,\di\rho_1\leq\int\phi\,\di\rho_2$ for all
  $\phi\in\Ci_{\rm cv}\big(\Pc(\ov S)\big)$.
\item There exists an $S$-valued random variable $X$ defined on some
  probability space $(\Om,\Fi,\P)$ and sub-\si-fields $\Fi_1\sub\Fi_2\sub\Fi$
  such that $\dis\rho_i=\P\big[\P[X\in\,\cdot\,|\Fi_i]\in\,\cdot\,\big]$
  $(i=1,2)$.
\end{enumerate}
\et
\bpro
For any probability kernel $P$ on $\Pc(\ov S)$, measure $\rho\in\Pc(\ov S)$,
and function $\phi\in\Ci(\Pc(\ov S))$, we define $\rho P\in\Pc(\Pc(\ov S))$ and
$P\phi\in B(\Pc(\ov S))$ by
\be
\rho P:=\int\rho(\di\mu)P(\mu,\,\cdot\,)
\quand
P\phi:=\int P(\,\cdot\,,\di\mu)\phi(\mu).
\ee
By definition, a \emph{dilation} is a probability kernel $P$ such that
$Pl_f=l_f$ for all $f\in\Ci(\ov S)$.
%The definition given by Strassen is that $P\phi=\phi$ for all continuous
%affine functions $\phi:\Ci(\ov S)\to\R$, but in the light of
%Lemma~\ref{L:conaf} we can reformulate this as done here.

As in the proof of Lemma~\ref{L:convrep}, we can view $\Pc(\ov S)$
as a convex compact metrizable subset of the locally convex topological
vector space $\Ci(\ov S)'$. Then \cite[Thm~2]{Str65} tells us that (i) is
equivalent to:
\begin{enumerate}\addtocounter{enumi}{2}
\item There exists a dilation $P$ on $\Pc(\ov S)$ such that $\rho_2=\rho_1P$.
\end{enumerate}
To see that this implies (ii), let $\xi_1,\xi_2$ be $\Pc(\ov S)$-valued random
variables such that $\xi_1$ has law $\rho_1$ and the conditional law of
$\xi_2$ given $\xi_1$ is given by $P$. Let $\Fi_1$ be the \si-field generated
by $\xi_1$, let $\Fi_2$ be the \si-field generated by $(\xi_1,\xi_2)$,
and let $X$ be an $\ov S$-valued random variable whose conditional law given
$\Fi_2$ is given by $\xi_2$. Then
\be
\P\big[\P[X\in\,\cdot\,|\Fi_2]\in\,\cdot\,\big]
=\P[\xi_2\in\,\cdot\,]=\rho_1P=\rho_2.
\ee
Since $P$ is a dilation
\be
\E[f(X)\,|\,\Fi_1]=\E\big[\E[f(X)\,|\,\Fi_2]\,\big|\,\Fi_1\big]
=\E\big[l_f(\xi_2)\,\big|\,\Fi_1\big]
=\int P(\xi_1,\di\mu)l_f(\mu)=l_f(\xi_1)
\ee
for all $f\in\Ci(\ov S)$, and hence
\be
\P\big[\P[X\in\,\cdot\,|\Fi_1]\in\,\cdot\,\big]
=\P\big[\xi_1\in\,\cdot\,\big]=\rho_1.
\ee
We note that since $\rho_1,\rho_2\in\Pc(\Pc(S))$, we have $\xi_1,\xi_2\in\Pc(S)$
a.s.\ and hence $X\in S$ a.s. This proves the implication (iii)$\volgt$(ii). 

To complete the proof, it suffices to show that (ii)$\volgt$(i). By
Lemma~\ref{L:convrep}, each $\phi\in\Ci_{\rm cv}(\Pc(\ov S))$ is of the form
$\phi=\sup_{f\in C}l_f$ for some $C\sub\Ci(\ov S)$. Then (ii) implies
\be\ba{l}
\dis\int\phi\,\di\rho_1
=\E\big[\sup_{f\in C}\E[f(X)\,|\,\Fi_1]\big]
=\E\big[\sup_{f\in C}\E\big[\E[f(X)\,|\,\Fi_2]\,\big|\,\Fi_1\big]\big]\\[5pt]
\dis\quad\leq
\E\big[\E\big[\sup_{f\in C}\E[f(X)\,|\,\Fi_2]\,\big|\,\Fi_1\big]\big]
=\E\big[\sup_{f\in C}\E[f(X)\,|\,\Fi_2]\big]
=\int\phi\,\di\rho_2.
\ec
\QED \epro

The $n$-th moment measure $\rho^{(n)}$ associated with a probability law
$\rho\in\Pc(\Pc(\ov S))$ has been defined in (\ref{momdef}). The following lemma
links the first and second moment measures to the convex order.

\bl[First and second moment measures]
Let\label{L:rho12} $S$ be a Polish space. Assume that
$\rho_1,\rho_2\in\Pc(\Pc(S))$ satisfy $\rho_1\leq_{\rm cv}\rho_2$. Then
$\rho^{(1)}_1=\rho^{(1)}_2$ and
\be\label{2comp}
\int\rho^{(2)}_1(\di x,\di y)f(x)f(y)\leq\int\rho^{(2)}_2(\di x,\di y)f(x)f(y)
\qquad\big(f\in B(S)\big).
\ee
If $\rho_1\leq_{\rm cv}\rho_2$ and (\ref{2comp}) holds with equality for all
bounded continuous $f:S\to\R$, then $\rho_1=\rho_2$.
\el
\bpro
By Theorem~\ref{T:Stras}, there exists an $\ov S$-valued random variable $X$
defined on some probability space $(\Om,\Fi,\P)$ and sub-\si-fields
$\Fi_1\sub\Fi_2\sub\Fi$ such that
$\dis\rho_i=\P\big[\P[X\in\,\cdot\,|\Fi_i]\in\,\cdot\,\big]$ $(i=1,2)$.
Since for each $f\in B(S)$ 
\be
\int\rho^{(1)}_1(\di x)f(x)=\E\big[\E[f(X)\,|\Fi_1]\big]=\E[f(X)]
=\E\big[\E[f(X)\,|\Fi_2]\big]=\int\rho^{(1)}_2(\di x)f(x),
\ee
we see that $\rho^{(1)}_1=\rho^{(1)}_2$. Fix $f\in B(S)$ and set
$M_i:=\E[f(X)\,|\Fi_i]$ $(i=1,2)$. Then
\be\ba{l}\label{M12}
\dis\int\rho^{(2)}_2(\di x,\di y)f(x)f(y)
=\E\big[\E[f(X)\,|\Fi_2]^2\big]
=\E[M_2^2]\\[5pt]
\dis\quad=\E[M_1^2]+\E\big[(M_2-M_1)^2\big]
\geq\E[M_1^2]=\int\rho^{(2)}_1(\di x,\di y)f(x)f(y),
%Here we have used that $\E[(M_2-M_1)M_1]
%=\E\big[\E[(M_2-M_1)\,|\,\Fi_1]M_1\big]=0.
\ec
proving (\ref{2comp}). Let $\ov S$ be a metrizable compactification of $S$.
If $\rho_1\leq_{\rm cv}\rho_2$ and (\ref{2comp}) holds with equality for all
bounded continuous $f:S\to\R$, then (\ref{M12}) tells us that $M_1=M_2$ for each
$f\in\Ci(\ov S)$, i.e.,
\be
\E[f(X)\,|\Fi_1]=\E[f(X)\,|\Fi_2]\mbox{ a.s.\ for each }f\in\Ci(\ov S).
\ee
By Lemma~\ref{L:separ}, we can choose a countable dense set
$\Di\sub\Ci(\ov S)$. Then $\E[f(X)\,|\Fi_1]=\E[f(X)\,|\Fi_2]$ for all 
$f\in\Di$ a.s.\ and hence $\P[X\in\,\cdot\,|\Fi_1]=\P[X\in\,\cdot\,\,|\Fi_2]$
a.s., proving that $\rho_1=\rho_2$.
\QED \epro

The following lemma shows that the convex order is a partial order,

\bl[Convex functions are distribution determining]
If\label{L:disdet} $\rho_1,\rho_2\in\Pc(\Pc(\ov S))$ satisfy
$\int\phi\,\di\rho_1=\int\phi\,\di\rho_2$ for all $\phi\in\Ci_{\rm
  cv}(\Pc(\ov S))$, then $\rho_1=\rho_2$.
\el
\bpro
For any $f\in\Ci(\ov S)$ and $\rho\in\Pc(\Pc(\ov S))$,
\be\ba{l}
\dis\int_{\ov S^2}\rho^{(2)}(\di x,\di y)f(x)f(y)\\[5pt]
\dis\quad=\int_{\Pc(\ov S)}\rho(\di\mu)\int_{\ov S^2}\mu(\di x)\mu(\di y)f(x)f(y)
=\int_{\Pc(\ov S)}\rho(\di\mu)l_f(\mu)^2.
\ec
Therefore, since $l_f^2$ is a convex function,
$\int\phi\,\di\rho_1=\int\phi\,\di\rho_2$ for all
$\phi\in\Ci_{\rm cv}(\Pc(\ov S))$ implies equality in (\ref{2comp}) and hence,
by Lemma~\ref{L:rho12}, $\rho_1=\rho_2$.
\QED \epro

\section{Open Problem~12 of Aldous and Bandyopadhyay}\label{A:problem}

We have seen that the use of the higer-level map from Section~\ref{S:hilev}
and properties of the convex order lead to an elegant and short proof of
Theorem~\ref{T:bivar}, which is similar to \cite[Thm~11]{AB05}. The most
significant improvement over \cite[Thm~11]{AB05} is that the implication
(ii)$\volgt$(i) is shown without a continuity assumption on the map $T$,
solving Open Problem~12 of \cite{AB05}. If one is only interested in solving
this open problem, taking the proof of \cite[Thm~11]{AB05} for granted, then
it is possible to give a shorter argument that does not involve the
higer-level map and the convex order.

One way to prove the implication (ii)$\volgt$(i) in Theorem~\ref{T:bivar} is
to show that nonendogeny implies the existence of a measure
$\nu\in\Pc(S^2)_\mu$ such that $T^{(2)}(\nu)=\nu$ and
$\nu\neq\ov\mu^{(2)}$. In \cite{AB05}, such a $\nu$ was constructed as the
weak limit of measures $\nu_n$ which satisfied $T^{(2)}(\nu_n)=\nu_{n+1}$;
however, to conclude that $T^{(2)}(\nu)=\nu$ they then needed to assume the
continuity of $T^{(2)}$. Their Open Problem~12 asks if this continuity
assumption can be removed.

In our proof of Theorem~\ref{T:bivar}, we take $\nu=\un\mu^{(2)}$, which by
Theorem~\ref{T:bivar2} and Lemma~\ref{L:rho12} from Appendix~\ref{A:cv}
satisfies $\nu\neq\ov\mu^{(2)}$ if and only if the RTP corresponding to $\mu$
is not endogenous, and by Lemma~\ref{momeas} satisfies $T^{(2)}(\nu)=\nu$.

Antar Bandyopadhyay told us that shortly after the publication of \cite{AB05},
he learned that their Open Problem~12 could be solved by adapting the proof of
the implication (3)$\volgt$(2) of \cite[Th\'eor\`eme~9]{BL07} to the setting
of RTPs. To the best of our knowledge, this observation has not been
published. The setting of \cite[Th\'eor\`eme~9]{BL07} are positive recurrent
Markov chains with countable state space, which are a very special case of the
RTPs we consider. In view of this, we sketch their argument here in our
general setting and show how it relates to our argument.

Let $(\om_\ibf,X_\ibf)_{\ibf\in\T}$ be an RTP corresponding to 
 the map $\gamma$ and a solution $\mu$
of a RDE. Construct $(Y_\ibf)_{\ibf\in\T}$ such that $(X_\ibf)_{\ibf\in\T}$
and $(Y_\ibf)_{\ibf\in\T}$ are conditionally independent and identically
distributed given $(\om_\ibf)_{\ibf\in\T}$. Then $X_\wurz=Y_\wurz$
a.s.\ if and only if the RTP corresponding to $\mu$ is endogenous.
Let $\nu$ denote the law of $(X_\wurz,Y_\wurz)$. Then
$\nu=\ov\mu^{(2)}$ if and only if endogeny holds. In view of this, to
  prove the implication (ii)$\volgt$(i) in Theorem~\ref{T:bivar}, it suffices
  to show that $\nu$ solves the bivariate RDE $T^{(2)}(\nu)=\nu$. This will
  follow provided we show that
\be
\big(\om_\ibf,(X_\ibf,Y_\ibf)\big)_{\ibf\in\T}
\ee
is an RTP corresponding to the map $\ga^{(2)}$ and $\nu$, i.e.,
\be\ba{rl}\label{RTP2}
{\rm(i)}&\mbox{the $(\om_\ibf)_{\ibf\in\T}$ 
are i.i.d.,}\\[5pt]
{\rm(ii)}&\mbox{for each $t\geq 1$, the $(X_\ibf,Y_\ibf)_{\ibf\in{\pa\T_{(t)}}}$ are
  i.i.d.\ with common law $\nu$}\\
&\mbox{and independent of $(\om_\ibf)_{\ibf\in\T_{(t)}}$,}\\[5pt]
{\rm(iii)}&\dis(X_\ibf,Y_\ibf)
=\ga^{(2)}[\om_\ibf]\big((X_{\ibf1},Y_{\ibf1}),\ldots,(X_{\ibf k_\ibf},Y_{\ibf \kappa(\om_\ibf)})\big)
\qquad(\ibf\in\T).
\ec
Here (i) and (iii) are trivial. To prove property~(ii), set
\be
\La_\ibf:=\big(X_\ibf,(\om_{\ibf\jbf})_{\jbf \in \T}\big)
\quand
\La^{(2)}_\ibf:=\big(X_\ibf,Y_\ibf,(\om_{\ibf\jbf})_{\jbf \in \T}\big)
\qquad(\ibf\in\T).
\ee
Then the $(\La_{\ibf})_{\ibf\in\T}$ are identically distributed. Moreover,
for each $t\geq 1$, the $(\La_{\ibf})_{\ibf\in\pa\T_{(t)}}$ are independent of
each other and of $(\om_\kbf)_{\kbf\in\T_{(t)}}$. Recall that
$(X_\ibf)_{\ibf\in\T}$ and $(Y_\ibf)_{\ibf\in\T}$ are conditionally
independent and identically distributed given
$(\om_\kbf)_{\kbf\in\T}$. Since the conditional law of $X_\ibf$ given
$(\om_\kbf)_{\kbf\in\T}$ only depends on
$(\om_{\ibf\jbf})_{\jbf\in\T}$, the same is true for $Y_\ibf$. Using
this, it is not hard to see that the $(\La^{(2)}_{\ibf})_{\ibf\in\T}$ are
identically distributed and for each $t\geq 1$, the
$(\La^{(2)}_{\ibf})_{\ibf\in\pa\T_{(t)}}$ are independent of each other and of
$(\om_\kbf)_{\kbf\in\T_{(t)}}$, and this in turn implies (ii).

In fact, since the law of $(X_\wurz,Y_\wurz)$ is the second moment
measure of the random measure $\xi_\wurz$ from Proposition~\ref{P:hlRTP},
the measure $\nu$ constructed here is the same as our measure
$\un\mu^{(2)}$. Thus, our argument and the one from \cite{BL07} are both based
on the same solution of the bivariate RDE.

\subsection*{Acknowledgements}

We thank Wolfgang L\"ohr and Jan Seidler for their help with
Lemmas~\ref{L:conaf} and \ref{L:convrep}, respectively.
We thank David Aldous, Antar Bandyopadhyay, and Christophe Leuridan for
answering our questions about their work.


\begin{thebibliography}{GKRV09}

\bibitem[AB05]{AB05}
D.J. Aldous and A.~Bandyopadhyay.
A survey of max-type recursive distributional equations.
\emph{Ann.\ Appl.\ Probab.}~15(2) (2005), 1047--1110.

\bibitem[Ant06]{Ant06}
A.~Bandyopadhyay.
A necessary and sufficient condition for the tail-triviality of a recursive
tree process.
\emph{Sankhya}~68(1) (2006), 1--23.

\bibitem[BL07]{BL07}
J.~Brossard and C.~Leuridan.
Cha\^ines de Markov consructives index\'ees par $\Z$.
\emph{Ann.\ Probab.}~35(2) (2007), 715--731.

\bibitem[Bou58]{Bou58}
N.~Bourbaki.
{\em \'El\'ements de Math\'ematique. VIII. Part. 1: Les Structures
  Fondamentales de l'Analyse. Livre III: Topologie G\'en\'erale. Chap. 9:
  Utilisation des Nombres R\'eels en Topologie G\'en\'erale. 2i\'eme \'ed.}
Actualit\'es Scientifiques et Industrielles~1045. Hermann \& Cie,
Paris, 1958.

\bibitem[Car00]{Car00}
N.L.~Carothers.
\emph{Real Analysis.}
Cambridge University Press, 2000.

\bibitem[CV77]{CV77}
C.~Castaing and M.~Valadier.
\emph{Convex Analysis and Measurable Multifunctions.}
Springer, Berlin, 1977.

\bibitem[Geo11]{Geo11}
H.-O.~Georgii. %Hans-Otto
\emph{Gibbs measures and phase transitions. 2nd ed.}
De Gruyter, Berlin, 2011.

\bibitem[Mos01]{Mos01}
E.~Mossel. %Elchanan Mossel
Reconstruction on trees: beating the second eigenvalue.
\emph{Ann.\ Appl.\ Probab.}~11(1) (2001), 285--300.

\bibitem[PW96]{PW96}
J.G.~Propp and  D. B.~Wilson. Exact sampling with coupled Markov chains
and applications to statistical mechanics.
\emph{Random Structures Algorithms}~9 (1996), 223--252.

\bibitem[Ros59]{Ros59} 
M.~Rosenblatt. 
Stationary processes as shifts of functions of independent
random variables. \emph{J.\ Math.\ Mech.}~8 (1959), 665--681.

\bibitem[Str65]{Str65}
V. Strassen.
The existence of probability measures with given marginals.
\emph{Ann.\ Math.\ Stat.}~36 (1965), 423--439.

\end{thebibliography}
\end{document}